\documentclass[sn-mathphys,Numbered]{sn-jnl}

\usepackage{graphicx}%
\usepackage{multirow}%
\usepackage{amsmath,amssymb,amsfonts}%
\usepackage{amsthm}%
\usepackage{mathrsfs}%
\usepackage[title]{appendix}%
\usepackage{xcolor}%
\usepackage{textcomp}%
\usepackage{manyfoot}%
\usepackage{booktabs}%
\usepackage{algorithm}%
\usepackage{algorithmicx}%
\usepackage{algpseudocode}%
\usepackage{listings}%

\theoremstyle{thmstyleone}%
\theoremstyle{thmstyletwo}%

\theoremstyle{thmstylethree}%
\newtheorem{thm}{Theorem}[section]
\newtheorem{cor}[thm]{Corollary}

\newtheorem{lem}[thm]{Lemma}

\newtheorem{ex}{Example}[section]
\newcommand{\be}{\begin{equation}}
\newcommand{\ee}{\end{equation}}
\newcommand{\ben}{\begin{enumerate}}
\newcommand{\een}{\end{enumerate}}
\newcommand{\beq}{\begin{eqnarray}}
\newcommand{\eeq}{\end{eqnarray}}
\newcommand{\beqn}{\begin{eqnarray*}}
\newcommand{\eeqn}{\end{eqnarray*}}
\begin{document}

\title[$GB\widetilde{W}$ metrics]{Generalized Berwald Projective Weyl ($GB\widetilde{W}$) metrics}
\author*[1]{\fnm{Nasrin} \sur{Sadeghzadeh}}\email{nsadeghzadeh@qom.ac.ir}
\affil*[1]{\orgdiv{Department of Mathematics}, \orgname{University of Qom}, \orgaddress{\street{Alghadir Bld.}, \city{Qom},
\country{Iran}}}
\maketitle
\thanks{This article has been published in \textit{Proceedings of the National Academy of Sciences, India Section A: Physical Sciences}, DOI:10.1007/s40010-024-00896-6.}

\abstract{This paper introduces a new quantity in Finsler geometry, called the generalized Berwald projective Weyl ($GB\widetilde{W}$) metric. The $C$-projective invariance of these metrics is demonstrated, and it is shown that they constitute a proper subset of the class of generalized Douglas ($GDW$) metrics. The paper also proves that all $GDW$ metrics with vanishing Landsberg curvature are of R-quadratic type. The class of $GDW$ metrics contains all Finsler metrics of scalar curvature, which provides an extension of the well-known Numata's theorem.}


\section{Introduction}\label{sec1}

Finsler spaces possess numerous projective invariants that are widely recognized, including Douglas curvature, Weyl curvature \cite{An}, and generalized Douglas-Weyl (GDW) curvature \cite{Ba}. Additionally, Akbar-Zadeh introduced H-curvature as a C-projective invariant \cite{Ak}, and a new C-projectively invariant quantity, known as $\widetilde{W}$-curvature or C-projective Weyl curvature, has been introduced to characterize Finsler metrics of constant flag curvature ($n>2$) \cite{NaTa}.\\
The equations $D_{j}{^i}_{kl}=0$, $W_{j}{^i}_{kl}=0$ and $D_{j}{^i}_{kl\mid m}{y^m}={T_{jkl}}{y^i}$ for some tensors $T_{jkl}$, are projectively invariant. The symbol $\mid$ denotes the horizontal covariant derivatives with respect to the Berwald connection.\\
Put simply, if a Finsler metric $F$ satisfies one of the aforementioned equations, then any Finsler metric that is projectively equivalent to $F$ must satisfy the same equations.\\

The study of GDW metrics is presented in \cite{Sa}, where it is demonstrated that all Finsler metrics of scalar curvature (Weyl metrics) belong to this class of metrics. Furthermore, B$\acute{a}$cs$\acute{o}$, in \cite{Ba}, shows that GDW metrics are closed under projective relations. This paper introduces a novel quantity in Finsler geometry, called the generalized Berwald projective Weyl ($GB\widetilde{W}$) curvature, which is a C-projective invariant. For a manifold $M$, let GB$\widetilde{W}(M)$ denote the class of all Finsler metrics that satisfy
\[
B_j{^i}_{kl}=\beta_j{^i}_{kl}+{b_{jkl}}{y^i},
\]
for some tensors $b_{jkl}$ and $\beta_{j}{^i}_{kl}$; where $\beta_{j}{^i}_{kl\mid m}{y^m}=0$.
This leads to a natural question: "How extensive is GB$\widetilde{W}(M)$, and what types of interesting metrics belong to this class?" \\

It is evident that all Berwald metrics are included in this class; however, Berwald metrics are not considered C-projective invariants. It is noteworthy that the class of Berwald metrics is a proper subset of the class of $GB\widetilde{W}$ metrics, which is illustrated in the following example.
\begin{ex}\cite{Sh2}
The standard Funk metric on the Euclidean unit ball is defined as follows
\[
F=\frac{\sqrt{|y|^2-(|x|^2|y|^2-<x,y>^2)}}{1-|x|^2}+\frac{<x,y>}{1-|x|^2}, \quad y\in T_xB^n\simeq R^n
\]
where $< , >$ and $\mid . \mid$ denote the Euclidean inner product and norm on $R^n$, respectively. $F$ is a Randers metric of constant flag curvature, and consequently $GB\widetilde{W}$ metric. It is easily seen that $\beta$ is a closed 1-form. By a simple calculation, it follows that $H= 0$ while $F$ is not Berwald metric.
\end{ex}
On the other hand, it has been proven that all Finsler metrics of dimension $n>2$ with constant flag curvature ($\widetilde{W}=0$) belong to this class of Finsler metrics. As an illustration, consider the following Finsler metric $F_\varepsilon$ (Example \ref{ex2}), which represents a particular family of Bryant's metrics expressed in a local coordinate system. This metric is locally projectively flat and has a constant flag curvature, thereby satisfying $\widetilde{W}=0$. Consequently, it is a $GB\widetilde{W}$ metric.
\begin{ex}\cite{Sh3} \label{ex2}
Let $\varepsilon$ be an arbitrary number with $\mid \varepsilon\mid<1$. Let
\[
F_\varepsilon:=\frac{\sqrt{\Psi\big(\frac{1}{2}(\sqrt{\Phi^2+(1-\varepsilon^2)|y|^4}+\Phi)\big)+(1-\varepsilon^2)<x,y>^2}+\sqrt{1-\varepsilon^2}<x,y>}{\Psi},
\]
where,
\[
\Phi:=\varepsilon |y|^2+(|x|^2|y|^2-<x,y>^2),\quad \quad \Psi:=1+2\varepsilon |x|^2+|x|^4.
\]
$F_\varepsilon=F_\varepsilon(x,y)$ is a Finsler metric on $R^n$. Note that if $\varepsilon=1$, then $F_1=\alpha_{+1}$ is the spherical metric on $R^n$.
\end{ex}
Next, we notice the metrics, which are not of constant curvature (i.e. $\widetilde{W}\neq 0$), are still of $GB\widetilde{W}$ type. Based on Lemma \ref{ScalarGBtildW}, it is clear that such a metric could not be of scalar curvature in dimension $n>2$. In other words, every Finsler metric of scalar curvature and non-constant curvature is not of $GB\widetilde{W}$ type while it is $GDW$ metric. This shows that the class of $GB\widetilde{W}$ metric is the proper subset of the class of $GDW$. The following example illuminates this point.
\begin{ex}\cite{Sh3}
Let
\[
F=\frac{\sqrt{\Big(|x|^2<a,y>-2<a,x><x,y>\Big)^2+|y|^2\Big(1-|a|^2|x|^4\Big)}}{1-|a|^2|x|^4}
\]
\[
-\frac{\Big(|x|^2<a,y>-2<a,x><x,y> \Big)}{1-|a|^2|x|^4},
\]
where, $a$ is a constant vector in $R^n$, which has very important properties. It is of scalar curvature and isotropic S-curvature; however, the flag curvature and the S-curvature are not constant. We have
\[
S=(n+1)<a,x>F,
\]
Then by lemma 9.1 in \cite{Sh3}, one has $E_{ij}=\frac{n+1}{2}<a,x>F_{ij}$. Then $H\neq 0$; it is not of $GB\widetilde{W}$ type while it is of $GDW$ metric.
\end{ex}
It will be demonstrated below that the class of $GDW$ metrics encompasses the class of $GB\widetilde{W}$ metrics. However, it is important to note that there exist numerous $GDW$ metrics that are not of the $GB\widetilde{W}$ type.
\begin{thm} \label{GDW}
Every $GB\widetilde{W}$ metric is $GDW$ metric.
\end{thm}
\begin{thm}\label{INVA}
The class of GB$\widetilde{W}(M)$ on Finsler manifold $(M,F)$ is closed under C-projective changes.
\end{thm}
The subsequent section establishes that the class of $GB\widetilde{W}$ metrics encompasses all Finsler metrics with a constant flag curvature ($n>2$). Additionally, an extension of the well-known Numata's theorem is proven.
\begin{thm} \label{Numata}
All non-Riemannian $GDW$ metrics with vanishing Landsberg curvature are of R-quadratic type ($n>2$).
\end{thm}
\section{Preliminaries}
A Finsler metric on a manifold $M$ is  a nonnegative function $F$ on $TM$ with the following properties
\begin{enumerate}
\item  $F$ is
 $C^{\infty}$ on $TM\setminus \{0\}$;
\item
$F(\lambda y) =\lambda F(y)$, $\forall \lambda >0$, $\ y\in TM$;
\item  For each $y\in T_xM$,
the following quadratic form ${\bf g}_y$ on $T_xM$ is positive definite,
\begin{equation}
{\bf g}_y(u, v):= {1\over 2} \Big [ F^2(y+ s u + tv ) \Big ] \mid_{s, t=0}, \quad \quad u, v\in T_xM.
\end{equation}
\end{enumerate}
At each point $x\in M$, $F_x:= F\mid_{T_xM}$, is an Euclidean norm, if and only if ${\bf g}_y$ is independent of $y\in T_xM\setminus\{0\}$.
To measure the non-Euclidean feature of $F_x$, define
${\bf C}_y: T_xM \times T_xM \times T_xM \to R$ by
\begin{equation}
{\bf C}_y(u, v, w):=
{1\over 2} {d \over dt} \Big [ {\bf g}_{y+tw} (u, v) \Big ] \mid_{t=0}, \ \ \ \ \ \ u, v, w\in T_xM.
\end{equation}
The family ${\bf C}:=\{{\bf C}_y\}_{y\in TM\setminus \{0\} }$ is called the {\it Cartan torsion}.
A curve $c(t)$ is called a {\it geodesic} if it satisfies
\begin{equation}
{d^2 c^i\over dt^2} + 2 G^i (\dot{c}(t))=0,
\end{equation}
where, $G^i(y)$ denotes local functions on $TM$ given by
\begin{equation}
G^i(y):= \frac{1}{4} g^{il}(y) \{ \frac{\partial^2 F^2}{\partial x^k \partial y^l} y^k - \frac{\partial F^2}{\partial x^l}\},\quad y\in T_xM.\label{Gi}
\end{equation}
$G^i$s called the associated spray to $(M,F)$. \\
$F$ is called a {\it Berwald metric} if
$G^i(y)$ are quadratic in $y\in T_xM$ for all $x\in M$. Define
\[
B_y:T_xM\otimes T_xM\otimes T_xM\rightarrow T_xM
\]
\[
B_y(u,v,w)=B_j{^i}_{kl}u^j v^k w^l \frac{\partial}{\partial x^i},
\]
where,
$
{B_j^i}_{kl}=\frac{\partial^3 G^i}{\partial y^j \partial y^k \partial y^l},
$
and
\[
E_y:T_xM\otimes T_xM \rightarrow R
\]
\[
E_y(u,v)=E_{jk} u^j v^k,
\]
where, $E_{jk}=\frac{1}{2}B_j{^m}_{km}$, $u=u^i\frac{\partial}{\partial x^i}$,
$v=v^i \frac{\partial}{\partial x^i}$ and $w=w^i\frac{\partial}{\partial x^i}$.
$B$ and $E$ are called the Berwald curvature and the mean Berwald curvature, respectively.
$F$ is called a Berwald metric and weakly Berwald (WB) metric if $B=0$ and $E=0$, respectively \cite{Sh2}.
By means of E-curvature, we can define $\bar{E}$-curvature as follows
\[
\bar{E}_y:T_xM\otimes T_xM\otimes T_xM \longrightarrow R
\]
\[
\bar{E}_y(u, v, w):=\bar{E}_{jkl}(y)u^iv^jw^k=E_{jk\mid l}u^jv^kw^l.
\]
It is worth noting that $\bar{E}_{ijk}$ is not completely symmetric with respect to all three indices. To define the $H$-curvature, we take the covariant derivative of $E$ along geodesics. Specifically, $H_{ij}=E_{ij\mid m}y^m$,
\[
H_y: T_xM\otimes T_xM \longrightarrow R
\]
\[
H_y(u, v):=H_{ij}u^i v^j
\]
Let
\[
D_j{^i}_{kl}=B_j{^i}_{kl}-\frac{1}{n+1}\frac{\partial^3}{\partial y^j \partial y^k \partial y^l}(\frac{\partial G^m}{\partial y^m}y^i).
\]
It is easy to verify that $D:=D_j{^i}_{kl} dx^j\otimes \frac{\partial}{\partial x^i}\otimes dx^k \otimes dx^l$ is a well-defined tensor on slit tangent bundle $TM_0$.
We call $D$ the Douglas tensor, which is a non-Riemannian projective invariant. For two Finsler metrics, $F$ and $\bar{F}$, with the Geodesic coefficients $G^i$ and $\bar{G}^i$, respectively. The diffeomorphism $\phi: F \rightarrow \bar{F}$ is called projective mapping if there exists a positive homogeneous scalar function $P(x, y)$ of degree one which is satisfying
\[
\bar{G^i}=G^i+P y^i,
\]
where, $P=P(x,y)$ is positively $y$-homogeneous of degree one which called projective factor \cite{DShen} and \cite{Sh2}.
A projective transformation with projective factor $P$ would be C-projective if $Q_{ij}= 0$; where,
\be\label{QCproj}
Q_{ij}=\frac{\partial Q_j}{\partial y^i}-\frac{\partial Q_i}{\partial y^j},
\ee
\[
Q_i=\frac{\partial P}{\partial x^i}-G^m_i \frac{\partial P}{\partial y^m}-P \frac{\partial P}{\partial y^i}.
\]
One could easily show that
\be\label{D2}
D_j{^i}_{kl}=B_j{^i}_{kl}-\frac{2}{n+1}\{E_{jk}\delta^i_l+E_{jl}\delta^i_k+E_{kl}\delta^i_j+E_{jkl}y^i\},
\ee
where $E_{jkl}=\frac{\partial E_{jk}}{\partial y^l}$.\\
In \cite{We}, Weyl introduces a projective invariant for Riemannian metrics. Then Douglas extends Weyl’s projective invariant to Finsler metrics \cite{Do}. Finsler metrics with vanishing projective Weyl curvature are called Weyl metrics or $W$ metrics. In \cite{Zs}, Szab$\acute{o}$ proves that Weyl metrics are exactly Finsler metrics of scalar flag curvature. In \cite{NaTa}, a new C-projective invariant quantity is defined, namely $\widetilde{W}$ curvature, as follows
\[
\widetilde{W}^i_k=K^i_k-\frac{1}{1-n^2}\{y^i\widetilde{K}_{0k}-\delta^i_k \widetilde{K}_{00}\},
\]
where, $\widetilde{K}_{jk}=nK_{jk}+K_{kj}+y^r\frac{\partial K_{kr}}{\partial y^j}$.\\
$\widetilde{W}{^i}_k$ is called projective Weyl curvature or $\widetilde{W}$ curvature which is another candidate for characterizing the Finsler metrics of constant flag curvature.\\
There is another projective invariant equation in Finsler geometry for some tensors $T_{jkl}$
\[
D_j{^i}_{kl\mid m}y^m=T_{jkl}y^i,
\]
where, $D_j{^i}_{kl\mid m}$ denotes the horizontal covariant derivatives of $D_j{^i}_{kl}$ with respect to the Berwald connection of $F$. These metrics
are called GDW metrics.\\
Now a new C-projective invariant equation is introduced as follows
\be \label{GBW}
B_j{^i}_{kl}=\beta_j{^i}_{kl}+b_{jkl}y^i,
\ee
for some tensors $b_{jkl}$ and $\beta_j{^i}_{kl}$ where $\beta_j{^i}_{kl\mid m}y^m=0$, or equivalently,
\[
h^i_rB_j{^r}_{kl\mid m}y^m=0.
\]
Finsler metrics satisfying \eqref{GBW} are called $GB\widetilde{W}$ metrics.\\
There are numerous Finsler metrics that belong to the GB$\widetilde{W}$ type. Specifically, all Berwald metrics and $\widetilde{W}$ metrics, as well as Finsler metrics with a constant flag curvature ($n>2$), are included in this class of metrics.
\section{Generalized Berwald Projective Weyl ($GB\widetilde{W}$) metrics}
This section focuses on $GB\widetilde{W}$ metrics. Firstly, it is demonstrated that all Berwald and $\widetilde{W}$ metrics belong to this noteworthy class of C-projective invariant metrics.
\lem \label{lem1}
The class of $GB\widetilde{W}$ metrics contains $\widetilde{W}$ metrics ($n>2$).
\proof
Let $F$ be a $\widetilde{W}$ metric ($n>2$). Based on \cite{NaTa}, it is of constant curvature. Then there is a constant $\lambda$ such that
\[
R^i_k=\lambda(F^2\delta^i_k-y^iy_k).
\]
Then
\[
R_j{^i}_{kl.m}=2\lambda (C_{jlm}\delta^i_k-C_{jkm}\delta^i_l).
\]
Here ${}_{.m}$ means the differential with respect to $y^m$. In \cite{Tayebi}, the following Ricci identity has been stated
\[
B_j{^i}_{ml\mid k}-B_j{^i}_{km\mid l}=R_{j}{^i}_{kl.m}.
\]
Contracting the above equations by $y^k$ one easily can find
\[
B_j{^i}_{ml\mid 0}=2\lambda C_{jlm}y^i.
\]
Thus, $h_m^i B_j{^m}_{kl\mid 0}=0$ which means that $F$ is of $GB\widetilde{W}$ type.
\endproof
\subsection{Proof of Theorem \ref{GDW}}
Let $F$ be a Finsler metric of $GB\widetilde{W}$ type; then the Berwald curvature satisfies the following equation
\[
B_j{^i}_{kl}=\beta_j{^i}_{kl}+b_{jkl}y^i,
\]
for some tensors $b_{jkl}$ and $\beta_{j}{^i}_{kl}$ where $\beta_j{^i}_{kl\mid 0}=0$. Then $B_j{^i}_{kl\mid 0}=b_{jkl\mid 0}y^i$ and contracting it by $y^l$ yields $b_{jkl\mid 0}y^l=0$.
Thus, by the definition of $E$-curvature, one gets
\[
2E_{jk}=\beta_j{^m}_{km}+b_{jkm}y^m.
\]
Thus, one may conclude that
\be
H_{jk}=E_{jk\mid 0}=0. \label{H}
\ee
Considering \eqref{D2}, one easily gets
\be\label{GD}
D_j{^i}_{kl}=(\beta_{j}{^i}_{kl}-\frac{2}{n+1}e_j{^i}_{kl})+(b_{jkl}+E_{jkl})y^i=t_j{^i}_{kl}+d_{jkl}y^i,
\ee
where,
\[
e_j{^i}_{kl}=E_{jk}\delta^i_l+E_{jl}\delta^i_k+E_{kl}\delta^i_j,\quad and \quad E_{jkl}=\frac{\partial E_{jk}}{\partial y^l}=\frac{1}{2}\frac{\partial^3 S}{\partial y^j \partial y^k \partial y^l}.
\]
According to \eqref{GD} and \eqref{H} and considering $\beta_j{^i}_{kl\mid m}y^m=0$, we will have $t_j{^i}_{kl\mid m}y^m=0$, meaning that $F$ is $GDW$ metric.
\endproof
\begin{cor}\label{1}
Let $F$ be a $GB\widetilde{W}$ metric, then $H=0$.
\end{cor}
\begin{cor}\label{2}
Every $GDW$ metric is $GB\widetilde{W}$ metric, if and only if $H=0$.
\end{cor}
\begin{cor}\label{3}
Every R-quadratic Finsler metric is $GB\widetilde{W}$ metric.
\end{cor}
The class of $\widetilde{W}$ metrics is a proper subset of the class of $GB\widetilde{W}$ metrics. Likewise, the class of $W$ metrics ($W=0$) is a proper subset of the class of $GDW$ metrics.
\\
\\
{\textbf{Proof of Theorem \ref{INVA}}
\\
Let $F$ be a Finsler metric of $GB\widetilde{W}$ type; then the Berwald curvature is as follows
\[
B_j{^i}_{kl}=\beta_j{^i}_{kl}+b_{jkl}y^i,
\]
for some tensors $b_{jkl}$ and $\beta_j{^i}_{kl}$ where $\beta_j{^i}_{kl\mid 0}=0$. Under a C-projective transformation, $\widetilde{G}^i=G^i+Py^i$, one has
\be
{\widetilde{B}}_{j}{^i}_{kl}=(\beta_j{^i}_{kl}+\rho_j{^i}_{kl})+(b_{jkl}+P_{jkl})y^i=\widetilde{\beta}_j{^i}_{kl}+\widetilde{b}_{jkl}y^i, \label{Ber}
\ee
where,
\[
\rho_{j}{^i}_{kl}=P_{jk}\delta^i_l+P_{jl}\delta^i_k+P_{kl}\delta^i_j,\quad P_{jk}=\frac{\partial^2 P}{\partial y^j \partial y^k}\quad and \quad P_{jkl}=\frac{\partial^3 P}{\partial y^j \partial y^k \partial y^l}.
\]
Based on the previous theorem and our understanding that the $GDW$ property is a projective invariant, we can conclude that both $F$ and $\widetilde{F}$ belong to the $GDW$ metric class. Therefore, taking into account equation \eqref{D2}, we obtain
\[
D_j{^i}_{kl}=(\beta_j{^i}_{kl}-\frac{2}{n+1}e_j{^i}_{kl})+(b_{jkl}-\frac{2}{n+1}E_{jkl})y^i=
\]
\[
\widetilde{D}_j{^i}_{kl}=(\widetilde{\beta}_j{^i}_{kl}-\frac{2}{n+1}\widetilde{e}_j{^i}_{kl})+(\widetilde{b}_{jkl}-\frac{2}{n+1}\widetilde{E}_{jkl})y^i,
\]
where,
\[
\widetilde{e}_j{^i}_{kl}=\widetilde{E}_{jk}\delta^i_l+\widetilde{E}_{jl}\delta^i_k+\widetilde{E}_{kl}\delta^i_j,\quad and \quad \widetilde{E}_{jkl}=\frac{\partial \widetilde{E}_{jk}}{\partial y^l}.
\]
Since $F$ belongs to both the classes of $GDW$ metric and $GB\widetilde{W}$ metric, we can conclude that
\[
\beta_j{^i}_{kl\mid 0}=\frac{2}{n+1}e_j{^i}_{kl\mid 0}=0,
\]
Moreover, $\widetilde{F}$ is $GDW$-metric; then
\be
\widetilde{\beta}_j{^i}_{kl\parallel 0}-\frac{2}{n+1}\widetilde{e}_j{^i}_{kl\parallel 0}=0, \label{FGDW}
\ee
where $\parallel$ denotes the horizontal derivative with respect to $\widetilde{G}^i$. Also, since it is shown $\widetilde{E}_{jk\parallel 0}=0$; then
\[
\widetilde{e}_j{^i}_{kl\parallel 0}=\widetilde{E}_{jk||0}\delta^i_l+\widetilde{E}_{jl\parallel 0}\delta^i_k+\widetilde{E}_{kl\parallel 0}\delta^i_j=0.
\]
Thus, based on equation \eqref{FGDW}, we can deduce that $\widetilde{\beta}_j{^i}_{kl||0}=0$, indicating that $\widetilde{F}$ is a $GB\widetilde{W}$ metric. Next, we aim to demonstrate that $\widetilde{E}_{jk\parallel 0}=0$. Using equation \eqref{Ber}, we obtain
\[
\widetilde{E}_{jk}=E_{jk}+\frac{n+1}{2}P_{jk},
\]
where, $P_{jk}=\frac{\partial^2 P}{\partial y^j \partial y^k}$. Then we will have
\[
\widetilde{E}_{jk\parallel 0}=y^m\partial_m \widetilde{E}_{jk}-2(G^m+Py^m)\widetilde{E}_{jk.m}-\widetilde{E}_{jm}(G^m_k+P_ky^m+P\delta^m_k)
\]
\[
-\widetilde{E}_{mk}(G^m_j+P_jy^m+P\delta^m_j)=y^m\partial_m \widetilde{E}_{jk}-2G^m\widetilde{E}_{jk.m}-\widetilde{E}_{jm}G^m_k-\widetilde{E}_{mk}G^m_j
\]
\be
=\widetilde{E}_{jk\mid 0}=E_{jk\mid 0}+\frac{n+1}{2}P_{jk\mid 0}.\label{EE}
\ee
Based on \eqref{QCproj}
\[
y^mQ_{jm.k}=y^m\Big[(\partial_j P_m-\partial_m P_j)-(G^r_jP_{rm}-G^r_mP_{rj})\Big]_{.k}
\]
\be
=-y^m\partial_m P_{jk}+2G^mP_{jk.m}+G^m_jP_{mk}+G^m_kP_{mj}=-P_{jk\mid 0}.\label{P0}
\ee
Based on the definition of C-projective transformation, which stated in \eqref{QCproj}, and considering \eqref{EE}, \eqref{P0} and \eqref{H}, one can conclude that $\widetilde{E}_{jk\parallel 0}=E_{jk\mid 0}=0$.
\endproof
\lem \label{ScalarGBtildW}
There does not exist a non-constant scalar flag curvature Finsler metric that is a $GB\widetilde{W}$ metric ($n>2$).
\proof
Assuming the existence of a $GB\widetilde{W}$ metric $F$ with scalar flag curvature, Corollary \ref{1} implies that the H-curvature would vanish. This would imply that $F$ has constant flag curvature, as demonstrated in \cite{Akbarzadeh}, which states that for a Finsler metric with scalar flag curvature in dimensions $n>2$, the flag curvature is constant if and only if $H=0$. However, this contradicts our initial assumption that $F$ has non-constant flag curvature, thus proving the non-existence of a $GB\widetilde{W}$ metric with scalar flag curvature that is non-constant.
\\
{\textbf{Proof of Theorem \ref{Numata}}
\\
Suppose $F$ is a $GDW$ metric with a Landsberg curvature of zero ($n>2$). In that case, it also has a vanishing $H$-curvature, as shown in \cite{Tayebi}. Utilizing Corollary \ref{2}, we can easily conclude that $F$ belongs to the $GB\widetilde{W}$ metric class.\\
Now assume that $F$ is of scalar curvature. Then by Lemma \ref{ScalarGBtildW}, it is of constant flag curvature $\lambda$. By Lemma \ref{lem1}, we have
\[
B_j{^i}_{kl\mid 0}=2\lambda C_{jkl}y^i.
\]
Contracting the above equation by $y_i$ yields $\lambda C_{jkl}=0$. $F$ is not Riemannian, then $\lambda=0$ and $R^i{_k}=0$. Then it is of R-quadratic type.\\
Now let $F$ be of non-scalar curvature. For this $GB\widetilde{W}$ metric $F$, we have
\[
B_j{^i}_{kl}=\beta_j{^i}_{kl}+b_{jkl} y^i,
\]
for some tensors $b_{jkl}$ and $\beta_j{^i}_{kl}$ where $\beta_j{^i}_{kl\mid 0}=0$. Since $F$ is of Landsberg type, then
\[
0=y_iB_j{^i}_{kl}=-2L_{jkl}=y_i\beta_j{^i}_{kl}+b_{jkl}F^2,
\]
Accordingly, one could get the following
\be\label{BerLand}
B_j{^i}_{kl}=\beta_j{^m}_{kl}h^i_m,
\ee
However, since $F$ is of $GB\widetilde{W}$ type, then $\beta_j{^i}_{kl\mid 0}=0$; then
\be \label{Bland}
B_j{^i}_{kl\mid 0}=0,
\ee
Taking the vertical derivative of the above equation with respect to $y^m$, we obtain the following
\be \label{Bland2}
B_j{^i}_{kl\mid m}+B_j{^i}_{kl\mid p.m}y^p=0.
\ee
But one has
\[
B_j{^i}_{kl\mid p.m}y^p=B_j{^i}_{kl.m;p}y^p-2G^rB_j{^i}_{kl.m.r}-N^r_mB_j{^i}_{kl.r}-N^r_jB_r{^i}_{kl.m}-N^r_kB_j{^i}_{rl.m}
\]
\[
-N^r_lB_j{^i}_{kr.m}+N^i_rB_j{^i}_{kl.m},
\]
where $B_j{^i}_{kl.m;p}=\frac{\partial}{\partial x^p} B_j{^i}_{kl.m}$. Then one can conclude that
\be\label{Bland3}
B_j{^i}_{kl\mid m.p}y^p=B_j{^i}_{km\mid l.p}y^p.
\ee
On the other hand, \eqref{Bland2} and \eqref{Bland3} yield
\be \label{Bland4}
B_j{^i}_{kl\mid m}=B_j{^i}_{km\mid l},
\ee
which putting in the following Ricci identity for Berwald connection \cite{Tayebi} yields
\[
R_j{^m}_{lk.p}=B_j{^m}_{pk\mid l}-B_j{^m}_{pl\mid k}=0.
\]
It means that $F$ is of R-quadratic type.
\endproof

\backmatter


{\small
{\em Authors' addresses}:
{\em Nasrin Sadeghzadeh}, University of Qom, Qom, Iran.\\ e-mail: \texttt{nsadeghzadeh@\allowbreak qom.ac.ir}.}

\end{document}